\documentclass[10pt,a4paper]{amsart}
\usepackage[T1]{fontenc}
\usepackage[english]{babel}
\usepackage{lmodern,amsmath,amssymb,mathtools,microtype}
\usepackage[a4paper,textwidth=170mm,textheight=229mm,centering]{geometry}
\usepackage{enumitem}
\setlist[enumerate]{label=\textup{(\roman*)},leftmargin=*,itemsep=1pt,topsep=3pt}
\newtheorem{thm}{Theorem}[section]
\newtheorem{Lemma}[thm]{Lemma}
\newtheorem{Proposition}[thm]{Proposition}
\newtheorem{Corollary}[thm]{Corollary}
\makeatletter
\let\c@subsection\c@thm
\makeatother

\usepackage{xcolor}
\usepackage{hyperref}

\newcommand{\eqnum}[1]{\refstepcounter{thm}\label{#1}\eqno(\thethm)}
\newcommand{\PP}{\mathbb P}
\newcommand{\QQ}{\mathbb Q}
\newcommand{\OO}{\mathcal O}
\DeclareMathOperator{\Pic}{Pic}
\DeclareMathOperator{\Gal}{Gal}
\DeclareMathOperator{\rk}{rk}
\DeclareMathOperator{\Tr}{Tr}
\DeclareMathOperator{\charac}{char}
\definecolor{DarkRed}{RGB}{173,0,0}
\definecolor{LightRed}{RGB}{201,0,0}

\hypersetup{
	colorlinks=true,
	linkcolor=DarkRed,
	citecolor=LightRed,
	urlcolor=LightRed
}

\title[On the unirationality of anticanonical weighted complete intersection surfaces]{On the unirationality of\\ anticanonical weighted complete intersection surfaces:\\ the infinite families}
\author[Giulia Gugiatti]{Giulia Gugiatti}
\address{Giulia Gugiatti\\ School of Mathematics, University of Edinburgh, Edinburgh EH9 3FD, United Kingdom}
\email{giulia.gugiatti@ed.ac.uk}
\author[Alex Massarenti]{Alex Massarenti}
\address{Alex Massarenti\\ Dipartimento di Matematica e Informatica, Universit\`a di Ferrara, Via Machiavelli 30, 44121 Ferrara, Italy}
\email{msslxa@unife.it}
\date{November 13, 2026}

\begin{document}
\begin{abstract}
Anticanonical  
weighted complete intersection surfaces fall into finitely many sporadic families and four infinite families. We prove that the anticanonical surfaces lying in the infinite families are unirational over their field of definition. Over a Hilbertian field, we further prove that a Hilbert subset of members consists of nonrational surfaces.
\end{abstract}

\maketitle
\tableofcontents

\section{Introduction}\label{sec:introduction}

Every log del Pezzo surface over a field of
characteristic zero is geometrically rational. However, over a field that is not
algebraically closed, geometric rationality does not imply
rationality. The exceptional curves on the minimal resolution may
form nontrivial Galois orbits, 
and a conic bundle may not admit a rational section
even when it has a rational multisection.
Log del
Pezzo surfaces that arise as weighted complete intersections have the
advantage of coming with explicit projective equations, which can be used in
birational constructions to establish rationality properties.

In this paper we work over a field $K$ of characteristic zero that is not algebraically closed and study anticanonical weighted complete intersection log del Pezzo surfaces over $K$ that belong to infinite families. These surfaces exhibit both phenomena described above. We prove
that every such surface is unirational over $K$. 
When $K$ is Hilbertian, we also prove that the members parametrised by a Hilbert subset of the parameter space are not rational over $K$.

We first recall the main definitions. 
A quasismooth \cite[Definition 6.9]{IF} and wellformed \cite[Definition 6.3]{IF}
weighted complete intersection
$X=X_{d_1,\ldots,d_c}\subset\PP_K(a_0,\ldots,a_n)$ is
\emph{anticanonical} if $\omega_X\simeq\OO_X(-1)$; by the adjunction formula \cite[Section 6.14]{IF}, this
means that $\sum_i a_i-\sum_j d_j=1$. We assume throughout that $X$ is not an
intersection with a linear cone, that is, $d_j\neq a_i$ for all $i,j$.
This excludes equations that eliminate a variable of the same weight,
and hence redundant presentations of the same surface.
Quasismoothness and wellformedness are geometric conditions: they hold
over $K$ if and only if they hold over its algebraic closure $\overline{K}$. The
classifications below, usually stated over $\mathbb{C}$, therefore
apply over $K$. We recall that wellformedness depends only on the degrees $d_i$ and the weights $a_j$, whereas quasismoothness depends on the defining equations \cite[Sections 6.10--12, Section 8]{IF}. By an \emph{infinite family} we mean a series of families $X_{d_1(k),\ldots,d_c(k)}\subset\PP_K(a_0(k),\ldots,a_n(k))$ 
indexed by an integer $k\geq 1$ such that, for every $k$, the family contains quasismooth and wellformed members.

Chen--Chen--Chen
\cite[Theorem~1.3]{ChenChenChen} show that every quasismooth log del
Pezzo weighted complete intersection that is not an intersection with
a linear cone has codimension $c\leq 2$. 
Johnson--Koll\'ar \cite[Theorem~8]{JohnsonKollar} classify 
anticanonical hypersurfaces: they fall into one infinite family and $22$
sporadic families. Kim--Park
\cite[Table~1 and Table~2]{KimPark} classify anticanonical codimension-two surfaces: they fall into three infinite families and
$38$ sporadic families. Together, these give the four infinite
families studied in this
paper:
$$
\begin{array}{c|c|c}
\text{abbreviation }&\text{weights}&\text{degrees}\\ \hline
\mathrm{JK}&(2,2k+1,2k+1,4k+1)&8k+4\\
\mathrm{KP}_1&(1,1,k+1,k+1,2k+1)&(2k+2,2k+2)\\
\mathrm{KP}_2&(1,2,2k+1,2k+1,4k+1)&(4k+2,4k+3)\\
\mathrm{KP}_3&(2,2k+1,2k+1,4k+1,6k+1)&(6k+3,8k+2)
\end{array}  \qquad k \in \mathbb{N}, k>1
\eqnum{eq:families}
$$ 
We stress that a family is determined by its weights and
degrees, and not every member is quasismooth and wellformed: the quasismooth and wellformed members
form a dense open subset of the parameter space. 

We often refer to each infinite family by the abbreviation in the first column. For fixed $k$, we write $\mathrm{JK}(k)$ and $\mathrm{KP}_i(k)$ for the
corresponding family.
The value $k=0$ is not included: for $\mathrm{JK}$ it gives  $X_4\subset\PP(2,1,1,1)$, whose quasismooth and wellformed members are smooth del Pezzo surfaces of degree two; for $\mathrm{KP}_1$ it gives $X_{2,2}\subset\PP^4$, which is listed separately among the sporadic codimension-two families; for $\mathrm{KP}_2$ and $\mathrm{KP}_3$ it gives intersections with a linear cone.

For a fixed row and integer $k$, let $B$ be the affine space of ordered tuples of weighted equations of the indicated degrees, and let $B^{\mathrm{qs}}$ be its quasismooth subset, i.e., the subset of tuples defining a
quasismooth surface. 
The \emph{generic member} is the surface over the function field $K(B)$ obtained by taking all coefficients independently.
If $K$ is Hilbertian, a member is \emph{Hilbert general} if its coefficients belong to a Hilbert subset of $B^{\mathrm{qs}}(K)$, as specified in Section~\ref{sec:arithmetic}. Such subsets are Zariski dense. 

Our main results are summarised in the following theorem.
\begin{thm}\label{thm:main}
Let $K$ be a field of characteristic zero.
Every anticanonical weighted complete intersection surface lying in an infinite family is $K$-unirational. More precisely:
\begin{enumerate}
\item Every quasismooth member of $\mathrm{JK}$ is birational to a degree-two del Pezzo surface with a $K$-rational generalized Eckardt point; this del Pezzo surface is the same for all $k \geq 1$. For each integer $k \geq 1$, every quasismooth member of $\mathrm{KP}_3(k)$ is birational to a quasismooth member of $\mathrm{JK}(k)$ and to a smooth cubic surface.
\item For each integer $k \geq 1$, every quasismooth member $X$ of $\mathrm{KP}_1(k)$ and of $\mathrm{KP}_2(k)$ has a smooth conic bundle model over $\PP^1_K$ with a rational bisection, reduced discriminant of degree $2k+5$, and canonical self-intersection $3-2k$.  Moreover, there is a dominant rational map $\PP^2_K\dashrightarrow X$ of degree two.
\item For each family and each integer $k\geq1$, the generic member is nonrational over $K(B)$. If $K$ is Hilbertian, a Hilbert subset of $B^{\mathrm{qs}}(K)$ consists of $K$-unirational surfaces that are not $K$-rational.
\end{enumerate}
\end{thm}

The unirationality assertion requires no hypothesis on rational points:
the equations themselves supply the necessary rational points and rational
curves. The nonrationality assertion has a different quantifier: it
concerns the generic member, defined over the field of independent
coefficients, and Hilbert specialisations. Every surface in the
theorem is geometrically rational.  
 
For the Johnson--Koll\'ar family, the substitution in Proposition~\ref{prop:JK-model} eliminates the parameter $k$: every
quasismooth member is birational over $K$ to a smooth del Pezzo surface
of degree two that does not depend on $k$. This agrees with the
geometric description in \cite{GugiattiRota}. 
The
resulting del Pezzo surface contains a distinguished rational point
lying on four geometric exceptional curves, that is,
a generalized Eckardt point. The unirationality criterion of
Salgado--Testa--V\'arilly-Alvarado
\cite[Corollary~3.3]{SalgadoTestaVarilly}, which refines Manin's
construction \cite[Theorem~29.4]{Manin}, excludes such points. Instead,
we use Gao's theorem \cite[Theorem~1]{Gao}, which treats exactly this
case and gives unirationality over $K$.

The first two Kim--Park families share a common geometric feature. Their equations determine a base-point-free pencil of binary quadratic forms. On a suitable $\PP^2$-bundle, the divisor obtained by setting one fiber coordinate equal to zero is the graph of the degree-two map defined by the pencil. It is a $K$-rational bisection of the conic bundle. Base change along this bisection makes the generic conic rational and gives the degree-two parametrisation in Theorem~\ref{thm:main}. For the second family, an initial model in $\PP^1\times\PP^2$ has a nonreduced fiber. We describe the $k+1$ elementary transformations that replace this fiber and compute the discriminant of the resulting smooth model. The first family is already smooth in the corresponding scroll. In both constructions, the discriminant has degree $2k+5$, which grows without bound as $k$ increases. For $k=1$ this is the degree-seven case covered by Koll\'ar--Mella \cite[Corollary~8]{KollarMella}; here the same explicit bisection works for every $k$. For comparison, unirational, nonrational conic bundles of discriminant degree eight over $C_1$ fields are constructed in \cite{CGM}.

Nonrationality is proved on smooth models. For the Johnson--Koll\'ar and third Kim--Park families, a quadratic twist together with a transitive permutation action on four exceptional curves, respectively
three lines, forces the generic Picard number to be one. The rationality criterion for minimal geometrically rational surfaces then applies. For the first two Kim--Park families, the universal discriminant and its cover parametrising the components of a singular fiber show that the generic conic bundle is standard. Its discriminant degree then gives nonrationality by Iskovskikh's conic bundle criterion. By Hilbert irreducibility, the relevant finite Galois actions and component covers are preserved under specialisation to
Hilbert general members. The arguments apply to the full coefficient spaces, not only to special subfamilies of equations.

\subsection*{Structure of the paper.} In Section~\ref{sec:normalforms} we give normal forms for the infinite families and describe the consequences of quasismoothness. In Section~\ref{sec:tools} we give the arithmetic criteria. In Sections~\ref{sec:JK}--\ref{sec:KP3} we construct the smooth models and prove unirationality. In Section~\ref{sec:arithmetic} we prove generic nonrationality and  the Hilbert specialisation statement. 

\subsection*{AI Use Statement.} The authors declare that no AI tools were used in the preparation of this paper.

\subsection*{Acknowledgments.}  GG acknowledges support from a Simons Investigator
Award (n.~929034). AM is a member of the Gruppo Nazionale per le Strutture Algebriche, Geometriche e le loro Applicazioni of the Istituto Nazionale di Alta Matematica ``F. Severi'' (GNSAGA-INdAM).

\section{Normal forms}\label{sec:normalforms}

We first place the infinite families in normal forms adapted to the birational constructions we use later. These forms isolate the terms relevant to quasismoothness and make explicit the resulting restrictions on the coefficients and geometric singularities.

\subsection*{The Johnson--Koll\'ar family}

Let $x,y_1,y_2,z$ be coordinates of weights $2,2k+1,2k+1,4k+1$.

\begin{Lemma}\label{lem:JK-normal}
Every quasismooth Johnson--Koll\'ar surface admits the equation
$$
a x^{4k+2}+b xz^2+x^{2k+1}q_2(y_1,y_2)+f_4(y_1,y_2)=0,
\qquad a,b\in K^*,
\eqnum{eq:JK-normal}
$$
where $q_2$ is binary quadratic and $f_4$ is a square-free binary quartic. It is well-formed. Its geometric singularities are $P_z$ of type $\frac1{4k+1}(1,1)$ and the four points of $(x=z=0)$ defined by the roots of $f_4$, each of type $\frac1{2k+1}(1,k)$.
\end{Lemma}
\begin{proof}
The weights give the general equation $a x^{4k+2}+b xz^2+x^{k+1}f_1(y_1,y_2)z+x^{2k+1}f_2(y_1,y_2)+f_4(y_1,y_2)=0$. If $a=0$ or $b=0$, the affine cone is singular at the point $(1,0,0,0)$ or $(0,0,0,1)$, respectively. A multiple root of $f_4$, including the case $f_4=0$, also gives a nonzero singular cone point with $x=z=0$. Thus $ab\ne0$ and $f_4$ is square-free. Setting $z'=z+x^kf_1/(2b)$ and dropping the prime gives \eqref{eq:JK-normal}, with $q_2=f_2-f_1^2/(4b)$.

The ambient space is well-formed, and its only positive-dimensional singular stratum is $(x=z=0)$. The surface meets it in four distinct points. At $P_z$, the equation solves for $x$ and leaves weights $(2k+1,2k+1)$ modulo $4k+1$. At a root of $f_4$, it solves for a coordinate along the singular line and leaves weights $(2,-1)$ modulo $2k+1$. These give the stated types; $P_x$ is not on the surface.
\end{proof}

\subsection*{The first Kim--Park family}

Set $m=k+1\geq2$ and let $x,y,z,t,w$ be coordinates of weights $1,1,m,m,2m-1$.

\begin{Lemma}\label{lem:KP1-normal}
Every quasismooth member of the first Kim--Park family has equations
$$
F_1=xw+P_1=0,\qquad F_2=yw+P_2=0,\qquad
P_i=Q_i(z,t)+A_{i1}(x,y)z+A_{i2}(x,y)t+H_i(x,y),
\eqnum{eq:KP1-normal}
$$
where $\deg Q_i=2$, $\deg A_{ij}=m$, and $\deg H_i=2m$. The quadratics $Q_1,Q_2$ have no common geometric root. The only geometric singularity is the $K$-rational point $P_w$, of type $\frac1{2m-1}(1,1)$.
\end{Lemma}
\begin{proof}
At the point $(0,0,0,0,1)$ of the affine cone, the two Jacobian rows have entries only in the $x,y$ columns. Quasismoothness makes the coefficient matrix of $xw,yw$ invertible. Row operations give \eqref{eq:KP1-normal}, and the weights give the remaining terms. If $Q_1,Q_2$ had a common root $[z_0:t_0]$, the two Jacobian rows at $(0,0,z_0,t_0,0)$ would be supported in the $z,t$ columns and would both annihilate $(z_0,t_0)$ by Euler's identity. Their rank would be at most one. The ambient singular line $(x=y=w=0)$ therefore misses the surface. At $P_w$, the equations solve for $x,y$ and leave weights $(m,m)$ modulo $2m-1$. Since $\gcd(m,2m-1)=1$, this is the stated quotient type. There are no other ambient singular strata meeting the surface.
\end{proof}

\subsection*{The second Kim--Park family}

Let $x,y,z,t,w$ be coordinates of weights $1,2,2k+1,2k+1,4k+1$.

\begin{Lemma}\label{lem:KP2-normal}
Every quasismooth member of the second Kim--Park family has equations \[F_1=xw+p=0 \quad  \text{and} \quad  F_2=yw+q=0,\] where
$$
\begin{aligned}
p&=Q_1(z,t)+xA_1(x,y)z+xA_2(x,y)t+H(x,y),\\
q&=xQ_2(z,t)+B_1(x,y)z+B_2(x,y)t+xR(x,y).
\end{aligned}
\eqnum{eq:KP2-normal}
$$
Here $Q_1,Q_2$ are binary quadratics, and $A_i,B_i,H,R$ have weighted degrees $2k,2k+2,4k+2,4k+2$, respectively. Moreover, $Q_1$ is square-free, $Q_1,Q_2$ have no common geometric root, and $H(0,1)\ne0$. The geometric singularities are $P_w$ of type $\frac1{4k+1}(1,1)$ and the two points of $(x=y=w=0)$ defined by the roots of $Q_1$, each of type $\frac1{2k+1}(1,k)$.
\end{Lemma}
\begin{proof}
The coefficient of $w$ in the first equation is a scalar multiple of $x$, with nonzero scalar by quasismoothness at $P_w$. Normalize it to $x$. The coefficient in the second is $\lambda x^2+\mu y$, with $\mu\ne0$ by the same Jacobian test. Subtracting $\lambda xF_1$ and rescaling gives $yw+q$. The expansion \eqref{eq:KP2-normal} follows from the weights.

At a common root of $Q_1,Q_2$ with $x=y=w=0$, the second Jacobian row vanishes. At a multiple root of $Q_1$, the first row vanishes. Both cases contradict quasismoothness. If $H(0,1)=0$, then $P_y$ lies on the surface and the first Jacobian row vanishes there. Thus $P_y$ is excluded. At $P_w$, the equations solve for $x,y$, leaving two coordinates of weight $2k+1$ modulo $4k+1$. At a root of $Q_1$, use a chart along the singular line: the two equations solve for the coordinate along the line and for $x$, since $Q_1$ is reduced and $Q_2$ does not vanish at that root. The remaining weights are $(2,-1)$ modulo $2k+1$. This proves the singularity statement.
\end{proof}

The pencil $[x^2:y]$ has base curve $\mathcal L=(x=y=Q_1(z,t)=0)\subset X$. Geometrically it consists of two weighted projective lines $\PP(2k+1,4k+1)$ meeting at $P_w$. The components and the two singular points on $(x=y=w=0)$ may be conjugate over $K$. In contrast, the pencil $[x:y]$ in the first family has base locus only $P_w$. No splitting of either quadratic pencil over $K$ will be assumed.

\subsection*{The third Kim--Park family}

Let $x,y_1,y_2,z,w$ be coordinates of weights $2,2k+1,2k+1,4k+1,6k+1$.

\begin{Lemma}\label{lem:KP3-normal}
Every quasismooth member of the third Kim--Park family admits equations
$$
xw+f_3(y_1,y_2)=0,\qquad
y_1w+bz^2+x^{2k}q_2(y_1,y_2)+a x^{4k+1}=0,
\qquad ab\ne0,
\eqnum{eq:KP3-normal}
$$
where $f_3$ is a square-free binary cubic with $f_3(0,1)\ne0$, and $q_2$ is binary quadratic.
\end{Lemma}
\begin{proof}
Quasismoothness at $P_w$ makes the coefficients of $w$ a nonzero multiple of $x$ and a nonzero linear form in $y_1,y_2$. Rescale the first equation and change the $y_i$ over $K$ to make these coefficients $x$ and $y_1$. The weights give $F_1=xw+f_3(y_1,y_2)+x^{2k+1}p_1(y_1,y_2)+c x^{k+1}z$ and $F_2=y_1w+bz^2+x^kq_1(y_1,y_2)z+x^{2k}q_2^0(y_1,y_2)+a x^{4k+1}$.
The coefficients $a,b$ are nonzero, since otherwise the second Jacobian row vanishes at $P_x$ or $P_z$. Set $w'=w+x^{2k}p_1+c x^kz$ and then $z'=z+x^k(q_1-cy_1)/(2b)$. These are invertible weighted changes. Dropping primes gives \eqref{eq:KP3-normal}, with $q_2=q_2^0-y_1p_1-(q_1-cy_1)^2/(4b)$. At a root of $f_3$ with $x=z=w=0$, the Jacobian rows are $df_3$ and $y_1\,dw$. Their independence forces $f_3$ to be square-free and excludes the root $[0:1]$.
\end{proof}

The same local calculation gives three geometric singularities of type $\frac1{2k+1}(1,k)$ at the points of $(x=z=w=0)$ defined by the roots of $f_3$ and the point $P_w$ of type $\frac1{6k+1}(1,3k+1)$. At $P_w$, the remaining weights after eliminating $x,y_1$ are $2k+1,4k+1$ modulo $6k+1$. The other ambient coordinate singularities are excluded by $a,b\ne0$. These observations also verify well-formedness for the normal forms above.

\section{Arithmetic tools}\label{sec:tools}

A del Pezzo surface is a smooth projective geometrically integral surface $S$ with $-K_S$ ample; its degree is $K_S^2$. In degree two, the anticanonical morphism $S\to\PP^2$ is a double cover, and $S$ is a smooth quartic in $\PP(2,1,1,1)$. A \emph{generalized Eckardt point} lies on four exceptional curves of $S_{\overline K}$; those curves need not be defined over $K$.

\begin{thm}\cite[Theorem~1]{Gao}\label{thm:Gao}
A degree-two del Pezzo surface over an infinite field of characteristic different from two is unirational if it has a rational generalized Eckardt point.
\end{thm}

For comparison, a nonconstant morphism $\PP^1_K\to S$ also gives unirationality of a degree-two del Pezzo surface in characteristic zero \cite[Theorem~3.2]{SalgadoTestaVarilly}. A rational point outside the ramification divisor that is not generalized Eckardt yields such a curve \cite[Corollary~3.3]{SalgadoTestaVarilly}. In general there is an extension of degree at most two over which $S$ is unirational \cite[Corollary~1.3]{SalgadoTestaVarilly}; Theorem~\ref{thm:Gao} lets us remain over $K$ for the models used here.

A \emph{conic bundle} in this paper is a flat projective morphism $\pi:S\to\PP^1_K$ from a smooth geometrically integral surface, with plane conics as fibers and smooth generic fiber. The bundle is \emph{standard} if $\rho_K(S/\PP^1)=1$. We use $K$-minimal for a smooth surface admitting no nontrivial birational morphism over $K$ to another smooth surface; this is distinct from relative minimality of a conic bundle.

\begin{Lemma}\label{lem:smooth-conic-discriminant}
Let a surface $S$ in a $\PP^2$-bundle over $\PP^1$ be given by a quadratic equation that is nonzero on each fiber. If $S$ is smooth, every singular geometric fiber has rank two and the determinant divisor is reduced. If this divisor has degree $\delta$, then $S$ is geometrically rational and $K_S^2=8-\delta$.
\end{Lemma}
\begin{proof}
The fibers have constant Hilbert polynomial, so the morphism is flat. A rank-one fiber would be a double line; all fiber derivatives vanish on that line, and the derivative in the base direction is a quadratic form on it. This form has a geometric zero, giving a singular point of $S$. Thus each singular fiber has rank two. At its vertex, smoothness of $S$ makes the derivative in the base direction nonzero. 

Diagonalizing the two nondegenerate directions locally then shows that the determinant has a simple zero. If the generic rank were at most two, every fiber would be singular. Over $\overline K$ the vertices $v(s)$ would then trace a curve on $S$; differentiating $q(s,v(s))=0$ and using $\partial_xq(s,v(s))=0$ would give $\partial_sq=0$ at the vertices, so $S$ would be singular there. Thus the generic fiber is smooth; flatness then gives geometric integrality of the total space. 

Over $\overline K$, Tsen's theorem gives a point on the generic conic \cite[Tag~03RD]{StacksTsen}, hence rationality. Contracting one component in each singular fiber yields a ruled surface over $\PP^1$, whose canonical square is $8$; each of the $\delta$ blow-ups reduces it by one.
\end{proof}

\begin{Proposition}[Enriques' criterion]\label{prop:Enriques}
If a conic bundle $S\to\PP^1_K$ has a rational multisection $\Gamma\simeq\PP^1_K$ of degree $e$, then there is a dominant rational map $\PP^2_K\dashrightarrow S$ of degree $e$.
\end{Proposition}
\begin{proof}
Base change by $\Gamma\to\PP^1$. The generic conic over $K(\Gamma)$ has the tautological point, so its function field is $K(\Gamma)(r)$. This is a rational function field in two variables over $K$. The original generic conic is geometrically integral, so its base change is integral and the projection to $S$ has degree $e$. A rational model of the base change gives the asserted map. In particular, a rational section makes $S$ rational.
\end{proof}

\begin{Proposition}\label{prop:standard-criterion}
For a conic bundle $S\to\PP^1_K$, the relative Picard number is one plus the number of closed discriminant points whose fibers split over their residue fields. Thus the bundle is standard exactly when every degenerate fiber is irreducible over the residue field of its base point.
\end{Proposition}
\begin{proof}
Over $\overline K$, let $F$ be a fiber class and choose a component $E_i$ from each of the $\delta$ singular fibers. Set $D_i=2E_i-F$. The classes $K_S,F,D_1,\ldots,D_\delta$ form a basis of $\Pic(S_{\overline K})_\QQ$: one has $D_i^2=-4$, $D_iD_j=0$ for $i\ne j$, and $D_iK_S=D_iF=0$, while the intersection form on $\QQ K_S+\QQ F$ has determinant $-4$. The geometric Picard rank is $\delta+2$ by the ruled-surface description in Lemma~\ref{lem:smooth-conic-discriminant}.

Galois acts on the $D_i$ by signed permutations. For an orbit of discriminant points, this representation has a one-dimensional invariant subspace precisely when the stabilizer of one point does not exchange the two components of its fiber. That condition says that the fiber splits over the residue field. Otherwise its invariant subspace is zero. The classes $K_S,F$ are invariant, and the rank of $\Pic(S)$ equals the dimension of the invariant geometric Picard group over $\QQ$. Subtracting the base Picard rank gives the formula.
\end{proof}

\begin{thm}[Iskovskikh's rationality criteria]\label{thm:nonrationality-criteria}
Over a field of characteristic zero, the following hold.
\begin{enumerate}
\item A smooth projective geometrically rational $K$-minimal surface is $K$-rational exactly when it has a $K$-point and $K_S^2\geq5$.
\item A standard conic bundle over $\PP^1_K$ with discriminant degree at least four is not $K$-rational. In particular, it has no rational section.
\end{enumerate}
\end{thm}
\begin{proof}
Part~\textup{(i)} is the minimal-surface criterion; see \cite[Section 4]{Iskovskikh1996} or \cite[Theorem~5.5]{Prokhorov}. For part~\textup{(ii)}, Lemma~\ref{lem:smooth-conic-discriminant} gives $K_S^2\leq4$. The case $K_S^2=4$ is \cite[Corollary~5.1.2(iii)]{Prokhorov}. For $K_S^2\leq3$, nonrationality is Iskovskikh's theorem \cite[Section 4]{Iskovskikh1996}; birational rigidity \cite[Corollary~5.1.2(i)]{Prokhorov} gives an alternative argument only when $K_S^2\leq0$. A rational section would contradict Proposition~\ref{prop:Enriques}.
\end{proof}

\section{The Johnson--Koll\'ar del Pezzo model}\label{sec:JK}

Following \cite{GugiattiRota}, we use the normal form of the Johnson--Koll\'ar family to construct a smooth degree-two del Pezzo model. We then identify on this model a distinguished $K$-rational generalized Eckardt point, which will yield unirationality.

\begin{Proposition}\label{prop:JK-model}
Let $X$ be a quasismooth Johnson--Koll\'ar surface in the normal form \eqref{eq:JK-normal}. Then $X$ is birational over $K$ to the smooth degree-two del Pezzo surface
$$
S_{\mathrm{JK}}:\quad G=a u^2+u\bigl(bw^2+q_2(v_1,v_2)\bigr)+f_4(v_1,v_2)=0
\subset\PP_K(2,1,1,1)_{u,v_1,v_2,w}.
\eqnum{eq:JK-dP2}
$$
The inverse maps on $x\ne0$ and $u\ne0$ are
$$
(u,v_1,v_2,w)=\left(x,\frac{y_1}{x^k},\frac{y_2}{x^k},\frac z{x^{2k}}\right),
\qquad
(x,y_1,y_2,z)=(u,u^kv_1,u^kv_2,u^{2k}w).
\eqnum{eq:JK-maps}
$$
The point $P=[0:0:0:1]\in S_{\mathrm{JK}}(K)$ is generalized Eckardt.
\end{Proposition}
\begin{proof}
The maps \eqref{eq:JK-maps} respect the weights and are inverse on the indicated opens. Substitution in \eqref{eq:JK-normal} gives $u^{4k}G$, so they also identify the localized affine cones. The cone over $S_{\mathrm{JK}}$ is therefore smooth on $u\ne0$. On $u=0$ with $(v_1,v_2)\ne(0,0)$, a derivative of the square-free form $f_4$ is nonzero. At the remaining nonzero cone points, $w\ne0$ and $G_u=bw^2\ne0$. Thus $S_{\mathrm{JK}}$ is quasismooth. It avoids the only singular point $[1:0:0:0]$ of the ambient space since $a\ne0$, and hence is smooth. Quasismoothness implies geometric integrality here: distinct positive-degree hypersurface components would meet outside the cone vertex and create a singular point. The corresponding dense opens identify $X$ and $S_{\mathrm{JK}}$ birationally.

Projection to $[v_1:v_2:w]$ is a finite flat double cover $\pi:S_{\mathrm{JK}}\to\PP^2$, since its equation is monic in $u$ after division by $a$. Weighted adjunction gives $-K_{S_{\mathrm{JK}}}=\pi^*\OO_{\PP^2}(1)$, an ample class of square two. Over $\overline K$, write $f_4=c\ell_1\ell_2\ell_3\ell_4$ with distinct linear forms. Each curve $E_i=(u=\ell_i(v_1,v_2)=0)$ maps isomorphically to a line under $\pi$. Thus $E_i\simeq\PP^1$, $-K_{S_{\mathrm{JK}}}\cdot E_i=1$, and adjunction gives $E_i^2=-1$. All four curves contain $P$, as required.
\end{proof}

\begin{Corollary}\label{cor:JK-unirational}
Every quasismooth Johnson--Koll\'ar surface over $K$ is $K$-unirational.
\end{Corollary}
\begin{proof}
Apply Theorem~\ref{thm:Gao} to the model and rational generalized Eckardt point supplied by Proposition~\ref{prop:JK-model}.
\end{proof}

The point $P$ lies outside the anticanonical ramification divisor, since $G_u(P)=b\ne0$, but it lies on four exceptional curves. Thus the good-point criterion recalled in Section~\ref{sec:tools} does not apply to this point; Gao's theorem does. The anticanonical involution exchanges $P$ with the second rational point $P'=[-b/a:0:0:1]$. No splitting assumption on $f_4$ is needed.

\section{Conic bundle models for the first two Kim--Park families}\label{sec:conics}

Set $m=k+1$ and $\mathcal T_m=\PP_{\PP^1}(\OO\oplus\OO\oplus\OO(m))$, using the convention fixed by Cox coordinates $a,b,v,z,t$ of degrees $\deg a=\deg b=(1,0)$, $\deg v=(-m,1)$, and $\deg z=\deg t=(0,1)$. The irrelevant ideal is $(a,b)\cap(v,z,t)$. Write $\xi$ for the class of $(z=0)$ and $f$ for the pullback of a point of the base. Then $(v=0)$ has class $\xi-mf$, and $\xi^3=m$, $\xi^2f=1$, $f^2=0$.

\subsection*{The first family}

Keep the notation of Lemma~\ref{lem:KP1-normal}. Eliminating $w$ gives $xP_2-yP_1=0$ in $\PP(1,1,m,m)$, birational to $X$ with inverse $w=-P_1/x$ over $x\ne0$. Set $C=aH_2(a,b)-bH_1(a,b)$ and $D_j=aA_{2j}(a,b)-bA_{1j}(a,b)$. Their degrees are $2m+1$ and $m+1$, respectively. The smooth model will be the divisor
$$
S_1:\quad Cv^2+D_1vz+D_2vt+aQ_2(z,t)-bQ_1(z,t)=0
\quad\text{in }\mathcal T_m.
\eqnum{eq:KP1-scroll}
$$
Its bidegree is $(1,2)$. On $a=v=1$, the equation becomes $P_2(1,b,z,t)-bP_1(1,b,z,t)=0$, and the inverse map to $X$ is $(b,z,t)\mapsto[1:b:z:t:-P_1(1,b,z,t)]$. Here $b=y/x$, and the fiber coordinates are the original $z/x^m,t/x^m$. The chart $b=v=1$ gives the corresponding identification over $y\ne0$. Thus $S_1\dashrightarrow X$ is birational and the projection $S_1\to\PP^1_{a,b}$ resolves the pencil $[x:y]$ birationally.

\begin{Proposition}\label{prop:KP1-conic}
For every quasismooth member of the first family, $S_1$ is smooth and its projection $\pi_1:S_1\to\PP^1$ is a conic bundle with reduced discriminant of degree $2k+5$.
\end{Proposition}
\begin{proof}
The charts with $v\ne0$ are smooth by their identification with $X\cap(x\ne0)$ and $X\cap(y\ne0)$. In \eqref{eq:KP1-scroll}, on $a=1$, a singular point with $v=0$ would satisfy $Q_2-bQ_1=0$ and $-Q_1=0$, from the equation and its base derivative. This contradicts the absence of a common root of $Q_1,Q_2$. The chart $b=1$ gives the same argument. No fiber equation is identically zero, since that would make $Q_1,Q_2$ proportional. Lemma~\ref{lem:smooth-conic-discriminant} now gives a conic bundle with reduced discriminant.

Write $aQ_2-bQ_1=\alpha z^2+\beta zt+\gamma t^2$, where $\alpha,\beta,\gamma$ are linear in $a,b$. Four times the determinant of the symmetric matrix in $(v,z,t)$ is
$$
\Delta_1=C(4\alpha\gamma-\beta^2)-D_1^2\gamma-D_2^2\alpha+D_1D_2\beta,
\qquad \deg\Delta_1=2m+3=2k+5.
\eqnum{eq:KP1-discriminant}
$$
The determinant is nonzero by Lemma~\ref{lem:smooth-conic-discriminant}. Every summand has the indicated homogeneous degree, which is therefore the degree of the discriminant divisor.
\end{proof}

\subsection*{The second family}

Keep the notation of Lemma~\ref{lem:KP2-normal}. Projection away from $w$ gives the birational hypersurface $Y=(xq-yp=0)\subset\PP(1,2,2k+1,2k+1)$, with inverse $w=-p/x$ over $x\ne0$. The pencil is $[a:b]=[x^2:y]$. On $a=1$, set $s=b/a$ and $y=sx^2$. Dividing the equation by $x^2$ gives
$$
c(s)u^2+d_1(s)uz+d_2(s)ut+Q_2(z,t)-sQ_1(z,t)=0,
\qquad u=x^{2k+1},
\eqnum{eq:KP2-affine-conic}
$$
where $c(s)=R(1,s)-sH(1,s)$ and $d_i(s)=B_i(1,s)-sA_i(1,s)$. The degrees are $\deg c=2m$ and $\deg d_i\leq m$; the leading coefficient of $c$ is $-H(0,1)\ne0$. Let $C$ be the homogenization of $c$ to degree $2m$, and let $\overline D_i$ be the homogenization of $d_i$ to degree $m$. A model in $\PP^1_{a,b}\times\PP^2_{u,z,t}$ is
$$
Z:\quad Cu^2+a^m\overline D_1uz+a^m\overline D_2ut
+a^{2m-1}(aQ_2-bQ_1)=0.
\eqnum{eq:KP2-product}
$$
The map from $Y$ is $[x:y:z:t]\mapsto([x^2:y],[x^{2k+1}:z:t])$. It is birational: on the dense chart $x\ne0$, normalize $x=1$ to recover $y,z,t$ from the two projective factors. The product model has a double-line fiber $u=0$ over $a=0$, since $C(0,1)=-H(0,1)\ne0$. It is singular along $a=u=0$. Its matrix determinant has degree $6m=6k+6$; in fact it equals $a^{4m-3}\Delta_2$, with $\Delta_2$ as in \eqref{eq:KP2-discriminant} below, so it vanishes to order $4m-2$ at $a=0$. This is not the discriminant degree of a smooth conic bundle.

The substitution $u=a^mv$ and removal of the common factor $a^{2m-1}$ in \eqref{eq:KP2-product} give
$$
S_2^{\mathrm{cb}}:\quad
aCv^2+a\overline D_1vz+a\overline D_2vt+aQ_2(z,t)-bQ_1(z,t)=0
\quad\text{in }\mathcal T_m.
\eqnum{eq:KP2-scroll}
$$
This is again a divisor of bidegree $(1,2)$ and is birational to $X$ over $K$.

The substitution is a composition of $m=k+1$ elementary transformations of the ambient projective bundle. More explicitly, set $\mathcal P_r=\PP_{\PP^1}(\OO\oplus\OO\oplus\OO(r))$, with fiber coordinate $u_r$ of degree $(-r,1)$. The transformation from $\mathcal P_r$ to $\mathcal P_{r+1}$ is centered at the line $a=u_r=0$ and is given locally by $u_r=a u_{r+1}$. Blowing up that line makes the strict transform of the old fiber $\PP^2$ have normal bundle $\OO_{\PP^2}(-1)$; contracting it to a point gives the new $\PP^2$-bundle. 

The exceptional divisor is $\mathbb F_1$, mapping to the new fiber as the blow-up of a point. For $0\leq r\leq m-1$, the successive proper-transform equations are $Cu_r^2+a^{m-r}\overline D_1u_rz+a^{m-r}\overline D_2u_rt+a^{2m-1-2r}(aQ_2-bQ_1)=0$.
The final substitution $u_{m-1}=av$, followed by division by $a$, gives \eqref{eq:KP2-scroll}. Thus the nonreduced fiber is replaced by the rank-two fiber $Q_1(z,t)=0$. The grading identity $(-r,1)=(1,0)+(-(r+1),1)$ also describes each transformation in Cox coordinates.

\begin{Proposition}\label{prop:KP2-conic}
For every quasismooth member of the second family, $S_2^{\mathrm{cb}}$ is smooth and its projection $\pi_2:S_2^{\mathrm{cb}}\to\PP^1$ is a conic bundle with reduced discriminant of degree $2k+5$. The fiber over $a=0$ is singular and reduced.
\end{Proposition}
\begin{proof}
On $a=1,v\ne0$, set $v=1$. Equation~\eqref{eq:KP2-scroll} is the equation of $Y$ on $x=1$, by \eqref{eq:KP2-affine-conic}; this is the smooth chart $X\cap(x\ne0)$. At $a=1,v=0$, a singular point would give $Q_1=Q_2=0$ from the equation and its base derivative, which is impossible. On $b=1,a=0$, the equation is $Q_1(z,t)=0$. Away from $[v:z:t]=[1:0:0]$, a fiber derivative is nonzero since $Q_1$ is square-free. At that vertex, the base derivative is $C(0,1)=-H(0,1)\ne0$. This proves smoothness everywhere. Again the quadratic pencil prevents an identically zero fiber, so Lemma~\ref{lem:smooth-conic-discriminant} applies.

Set $aQ_2-bQ_1=\alpha z^2+\beta zt+\gamma t^2$ as before. Four times the determinant is
$$
\Delta_2=a\Bigl(C(4\alpha\gamma-\beta^2)
+a(-\overline D_1^2\gamma-\overline D_2^2\alpha+\overline D_1\overline D_2\beta)\Bigr),
\qquad\deg\Delta_2=2m+3.
\eqnum{eq:KP2-discriminant}
$$
The residual factor has degree $2m+2=2k+4$ and is nonzero at $a=0$, since $C(0,1)\ne0$ and $Q_1$ is square-free. Thus the fixed factor $a$ is simple, and the total discriminant degree is $2k+5$.
\end{proof}

\subsection*{Bisections, parametrisations, and nonrationality tests}

We now exploit the common geometry of the first two Kim--Park families. The conic bundle models carry a natural rational bisection, yielding explicit degree-two parametrisations, while their discriminants provide a convenient criterion for nonrationality over $K$.

\begin{Proposition}\label{prop:KP-bisections}
For either conic bundle in Propositions~\ref{prop:KP1-conic} and~\ref{prop:KP2-conic}, the curve $\Gamma=(v=0)\cap S$ is a rational bisection. Moreover, $K_S^2=3-2k$ and $\Gamma^2=-(2k+1)$. Every quasismooth member of either family admits a dominant rational map $\PP^2_K\dashrightarrow X$ of degree two.
\end{Proposition}
\begin{proof}
On $(v=0)\simeq\PP^1_{a,b}\times\PP^1_{z,t}$, both equations give $aQ_2-bQ_1=0$. Since $Q_1,Q_2$ have no common geometric root, this is the graph of the morphism $[z:t]\mapsto[Q_1(z,t):Q_2(z,t)]$. The second projection identifies $\Gamma$ with $\PP^1_K$, and the first has degree two. Proposition~\ref{prop:Enriques} gives the required parametrisation. Finally, both surfaces have class $2\xi+f$ in $\mathcal T_m$. Adjunction gives $K_S=(-\xi+(m-1)f)|_S$, hence $K_S^2=5-2m=3-2k$ and $\Gamma^2=(\xi-mf)^2(2\xi+f)=1-2m=-(2k+1)$. This agrees with the discriminant degrees in \eqref{eq:KP1-discriminant} and \eqref{eq:KP2-discriminant}.
\end{proof}

\begin{Corollary}\label{cor:KP-individual-nonrational}
For either of the first two Kim--Park families, if every singular fiber of the associated conic bundle is irreducible over the residue field of its base point, then $X$ is not $K$-rational. The conic bundle is standard and has no rational section.
\end{Corollary}
\begin{proof}
Proposition~\ref{prop:standard-criterion} gives standardness, and Theorem~\ref{thm:nonrationality-criteria}\textup{(ii)} applies since $\deg\Delta=2k+5\geq7$.
\end{proof}

In the second family, the criterion requires $Q_1$ to be irreducible over $K$ and each residual rank-two fiber in \eqref{eq:KP2-discriminant} to be nonsplit over its residue field. Even a fiber over a $K$-point may have two conjugate components. In particular, the bundle cannot be standard when $K$ is quadratically closed, since the fiber over $a=0$ then splits. For a standard member of either family, the generic conic has a nonzero Brauer class: a zero class would give a section and hence rationality. Section~\ref{sec:arithmetic} proves that the standardness criterion holds for the generic members and for Hilbert general specialisations.

\section{The third Kim--Park family}\label{sec:KP3}

The third Kim--Park family is birationally different from the first two. Its normal form yields both a Johnson--Koll\'ar model and a smooth cubic surface model. The former gives unirationality through the construction of Section~\ref{sec:JK}, while the latter will be used later to study nonrationality.

\begin{Proposition}\label{prop:KP3-models}
Let $X$ be a quasismooth member of the third family in the normal form \eqref{eq:KP3-normal}. Projection away from $w$ is a birational map to the quasismooth Johnson--Koll\'ar surface
$$
Y:\quad a x^{4k+2}+bxz^2+x^{2k+1}q_2(y_1,y_2)-y_1f_3(y_1,y_2)=0
\subset\PP_K(2,2k+1,2k+1,4k+1).
\eqnum{eq:KP3-JK}
$$
The surface $X$ is also birational to the smooth cubic
$$
S_3:\quad aUT^2+T\bigl(bW^2+q_2(U,V)\bigr)-f_3(U,V)=0
\subset\PP^3_{K,[U:V:W:T]}.
\eqnum{eq:KP3-cubic}
$$
On $xy_1\ne0$ and on $UT\ne0$, respectively, the latter map and its inverse are
$$
\begin{aligned}
[x:y_1:y_2:z:w]&\longmapsto
\left[\frac{y_1}{x^k}:\frac{y_2}{x^k}:\frac z{x^{2k}}:\frac{x^{k+1}}{y_1}\right],\\
[U:V:W:T]&\longmapsto
[\eta:\eta^kU:\eta^kV:\eta^{2k}W:-\eta^{3k-1}f_3(U,V)],\qquad \eta=UT.
\end{aligned}
\eqnum{eq:KP3-maps}
$$
\end{Proposition}
\begin{proof}
Eliminating $w$ by $xF_2-y_1F_1$ gives \eqref{eq:KP3-JK}, and $w=-f_3/x$ recovers the original surface on $x\ne0$. This also identifies the localized affine cones, so the cone over $Y$ is smooth there. On $x=0$ with $(y_1,y_2)\ne(0,0)$, the binary quartic $-y_1f_3$ is square-free by Lemma~\ref{lem:KP3-normal}, and one of its derivatives is nonzero. At the remaining nonzero cone points, $z\ne0$ and the derivative in $x$ is $bz^2\ne0$. Thus $Y$ is quasismooth.

By Proposition~\ref{prop:JK-model}, its smooth degree-two model is $a u^2+u(bW^2+q_2(U,V))-Uf_3(U,V)=0$. Set $u=UT$ and divide by $U$ on $U\ne0$ to obtain \eqref{eq:KP3-cubic}. The maps \eqref{eq:KP3-maps} follow and are inverse on the stated dense opens. The cubic is smooth on $U\ne0$ by its identification with the degree-two model. If a singular point had $U=0$, its equation and its $T$ derivative would give $f_3(0,V)=0$ and $bW^2+q_2(0,V)=0$. Since $f_3(0,1)\ne0$ and $b\ne0$, this forces $V=W=0$. Then $T\ne0$ and the $U$ derivative is $aT^2\ne0$, a contradiction.
\end{proof}

\begin{Corollary}\label{cor:KP3-unirational}
Every quasismooth member of the third Kim--Park family over $K$ is $K$-unirational.
\end{Corollary}
\begin{proof}
Apply Corollary~\ref{cor:JK-unirational} to the surface \eqref{eq:KP3-JK}. Thus this unirationality proof also passes through the degree-two model and Theorem~\ref{thm:Gao}.
\end{proof}

The cubic also has the rational points $[0:0:0:1]$ and $[0:0:1:0]$, so the cubic unirationality theorem \cite{Kol02} gives another proof. The latter point is Eckardt: its tangent plane is $T=0$, whose section is the union of the three lines above the roots of $f_3$. This cubic model, rather than the degree-two model, will be used to prove Picard number one. The degree-two model contains the $K$-defined exceptional curve $u=U=0$ and therefore has Picard number at least two. In fact, $[u:U:V:W]\mapsto[U^2:UV:UW:u]$ is a birational morphism from the degree-two model onto $S_3$. It contracts the $K$-defined exceptional curve $\{U=0,\ au+bW^2+q_2(0,V)=0\}$ to $[0:0:0:1]$, so the degree-two model is the blow-up of $S_3$ at this point.

\section{Generic nonrationality and Hilbert specialisations}\label{sec:arithmetic}

We use the coefficient spaces $B$ defined in Section~\ref{sec:introduction}. A field is \emph{Hilbertian} if irreducible polynomials over $K(t)$ that are separable in the polynomial variable remain irreducible for a Zariski dense set of specialisations $t\in K$. For a nonempty open subset of an affine coefficient space, a \emph{Hilbert subset} is obtained by imposing finitely many such irreducibility conditions, together with Zariski open conditions. Over a Hilbertian field it is Zariski dense; see \cite[Chapters~12--13]{FriedJarden}. Number fields are Hilbertian, as are the fields $F(t_1,\ldots,t_n)$ with $n\geq1$; in our characteristic-zero setting this includes $\mathbb C(t)$ and $\mathbb R(t)$.

For a smooth geometrically rational surface, $\rho_K(S)=\rk\Pic(S)$ equals the dimension of $\Pic(S_{\overline K})_\QQ^{\Gal(\overline K/K)}$. A del Pezzo surface of degree two or three with $\rho_K(S)=1$ is $K$-minimal, and hence nonrational by Theorem~\ref{thm:nonrationality-criteria}\textup{(i)}. We first establish this condition for the two del Pezzo models.

\begin{Lemma}\label{lem:twist-Picard}
Let $L$ be a field of characteristic zero and $S/L$ a smooth del Pezzo surface of degree $d\in\{2,3\}$ with an involution $\tau$ defined over $L$. Suppose that its geometric fixed locus consists of a smooth genus-one curve and $4-d$ isolated points. Assume that $\tau$ preserves $6-d$ exceptional curves with linearly independent classes, and that Galois acts transitively on these curves. If $b$ is transcendental over $L$, the twist of $S_{L(b)}$ by $\tau$ and $L(b)(\sqrt b)/L(b)$ has Picard number one over $L(b)$.
\end{Lemma}
\begin{proof}
The geometric Picard rank is $10-d$. The Lefschetz fixed-point formula gives $2+\Tr(\tau^*\mid\Pic(S_{\overline L})_\QQ)=4-d$; see \cite[Section 2.2]{DolgachevMartin}. Its invariant subspace therefore has dimension $6-d$, and the stated exceptional classes form a basis of that subspace.

Over $\overline L(\sqrt b)$ the twisted surface is identified with the constant surface $S_{\overline L}$. The automorphism $\sqrt b\mapsto-\sqrt b$ acts on its Picard group as $\tau$. Since $b$ is transcendental, this action and the automorphisms of $\overline L/L$ are independent. Consequently, a class defined over $L(b)$ must lie in the span of the exceptional classes and be invariant under their transitive permutation action. That invariant space is one-dimensional and contains the canonical class, proving the assertion.
\end{proof}

\begin{Proposition}\label{prop:JK-KP3-generic}
For each $k\geq1$, the generic Johnson--Koll\'ar member and the generic member of the third Kim--Park family are nonrational over their full coefficient function fields. If $K$ is Hilbertian, a Hilbert general quasismooth member of either family is not $K$-rational.
\end{Proposition}
\begin{proof}
For the Johnson--Koll\'ar model \eqref{eq:JK-dP2}, take $L$ to be the rational function field over $K$ in $a$ and the coefficients of $q_2,f_4$, with $b$ a further independent parameter. The change $w'=\sqrt b\,w$ identifies the surface with the twist of the $b=1$ model by $\tau:w'\mapsto-w'$. On that model, the fixed locus is the smooth anticanonical genus-one curve $w'=0$ and the two points $[0:0:0:1]$, $[-1/a:0:0:1]$. This curve is smooth, being the one-dimensional part of the fixed locus of an involution of a smooth surface; it has genus one since it lies in $|-K|$.

The four curves $E_i=(u=\ell_i(v_1,v_2)=0)$ above the roots of $f_4$ are individually preserved by $\tau$. They meet pairwise transversely at $P$ and nowhere else. Their intersection matrix is $J_4-2I_4$, of determinant $-16$, where $J_4$ has every entry equal to one. 

Their classes are independent. The generic binary quartic is irreducible over $L$, so Galois acts transitively on them. Lemma~\ref{lem:twist-Picard} gives Picard number one.

For the cubic \eqref{eq:KP3-cubic}, take $L=K(a,\operatorname{coeff}(q_2),\operatorname{coeff}(f_3))$ and again keep $b$ independent. Its $b=1$ model has involution $W\mapsto-W$, with fixed locus the smooth plane cubic $W=0$ and the point $[0:0:1:0]$. The three lines $T=\ell_i(U,V)=0$ above the roots of $f_3$ meet pairwise once at that point. Their intersection matrix is $J_3-2I_3$, of determinant $4$. They are individually preserved by the involution and form one Galois orbit, since the generic binary cubic is irreducible. The lemma again gives Picard number one. Theorem~\ref{thm:nonrationality-criteria}\textup{(i)} proves nonrationality of both smooth models.

These computations apply to the full generic weighted families, since the normalizations are rationally invertible after retaining the eliminated coefficients as free parameters. In the Johnson--Koll\'ar case, the inverse change is $f_2=q_2+f_1^2/(4b)$, with $f_1$ free. In the third Kim--Park case, retain the nonzero coefficients used to normalize $xw,y_1w$, together with $p_1,q_1,c$ from the proof of Lemma~\ref{lem:KP3-normal}; the inverse quadratic change is $q_2^0=q_2+y_1p_1+(q_1-cy_1)^2/(4b)$. The linear change in $y_1,y_2$ is rationally invertible on a coefficient chart. Thus the original coefficient fields are purely transcendental extensions of the normalized fields, and the geometric Picard action has the same invariant rank.

For specialisation, restrict the full coefficient space to a nonempty open where the models are smooth and the displayed changes are defined. The relative schemes of exceptional curves are finite \mbox{\'etale}, with $56$ and $27$ geometric points, respectively. 

Their classes generate the geometric Picard groups. After shrinking, the splitting field of all exceptional curves of the generic model defines a connected finite Galois \mbox{\'etale} cover. Hilbert irreducibility preserves its full Galois group on a Hilbert subset \cite[Chapter~12]{FriedJarden}. 

The marking identifies the intersection lattices and their Galois actions, so these specialisations retain Picard number one. The same rationality criterion applies, and the birational maps transfer nonrationality to the weighted surfaces.
\end{proof}

\begin{Proposition}\label{prop:KP-generic-standard}
For either of the first two Kim--Park families and every $k\geq1$, the conic bundle of the generic member over the full coefficient field is standard. If $K$ is Hilbertian, the same holds for a Hilbert general quasismooth member over $K$.
\end{Proposition}
\begin{proof}
Set $m=k+1$. First consider the effective coefficient spaces of the following two equations in $\mathcal T_m$:
$$
\begin{array}{ll}
\mathrm{I}:&Cv^2+D_1vz+D_2vt+aQ_2-bQ_1=0,
\quad\deg C=2m+1,\quad\deg D_i=m+1;\\
\mathrm{II}:&aCv^2+a\overline D_1vz+a\overline D_2vt+aQ_2-bQ_1=0,
\quad\deg C=2m,\quad\deg \overline D_i=m.
\end{array}
\eqnum{eq:universal-conics}
$$
All coefficients indicated here are independent. In either affine coefficient space $A$, restrict to the open where $Q_1,Q_2$ have no common root and the determinant divisor is reduced. A reduced determinant forces every singular fiber to have rank two with a simple zero of the determinant, and local diagonalization then shows that the total space is smooth; thus each member of this open is a conic bundle. This open is nonempty. In type~I, take $D_1=D_2=0$, $C=a^{2m+1}-b^{2m+1}$, $Q_1=z^2$, and $Q_2=t^2$. The discriminant is a nonzero multiple of $ab(a^{2m+1}-b^{2m+1})$. In type~II, take $\overline D_1=\overline D_2=0$, $C=a^{2m}-b^{2m}$, $Q_1=z^2+t^2$, and $Q_2=2z^2+3t^2$. The discriminant is a nonzero multiple of $a(a^{2m}-b^{2m})(2a-b)(3a-b)$. Both are square-free in characteristic zero. In type~II, reducedness at $a=0$ already forces $C(0,1)\ne0$ and $Q_1$ square-free, since the residual factor at $a=0$ is a nonzero multiple of $C(0,1)\operatorname{disc}(Q_1)$.

Let $U=\PP^1$ in type~I and $U=(a\ne0)$ in type~II. Before restricting the coefficient space, evaluation at a point of $U$ is surjective onto the six-dimensional space of ternary quadratic forms: each of its six coefficients can be varied independently. Locally on $U$, this is a surjective vector-bundle map. The universal incidence of rank-two quadrics is therefore an affine bundle over the rank-two locus of the determinant hypersurface, fibered over $U$. This locus is geometrically integral. Its double cover choosing one of the two component lines is geometrically integral as well: it is parametrised by an ordered pair of distinct lines and a nonzero scalar. Pullback through the evaluation map, and then restriction to the nonempty coefficient open, preserve geometric integrality of both incidences.

Over $K(A)$, the discriminant in $U$ is consequently one closed point, and its component cover is a nontrivial quadratic field extension. The degree of this discriminant point over $K(A)$ is $2m+3$ in type~I and $2m+2$ in type~II; the component cover has twice that degree. Type~II has the additional fixed fiber $Q_1(z,t)=0$ at $a=0$; its components are conjugate, since the discriminant of a generic binary quadratic is not a square. Thus every closed singular fiber is nonsplit over its residue field. Proposition~\ref{prop:standard-criterion} proves standardness.

For Hilbert specialisation, shrink $A$ so that the discriminant and its component cover are finite \mbox{\'etale} over $A$. In type~II, treat the fixed and residual discriminant components separately. For each integral discriminant component $\mathcal D$, its component cover $\widetilde{\mathcal D}$ is integral and has function field of degree $2[K(\mathcal D):K(A)]$ over $K(A)$. Hilbert irreducibility preserves these degrees and irreducibility simultaneously on a Hilbert subset. Each specialised discriminant component is then a single closed point, and its two fiber components remain conjugate over its residue field. This gives standardness after specialisation.

It remains to pass from \eqref{eq:universal-conics} to the full weighted coefficient spaces. In type~I, the linear maps $(H_1,H_2)\mapsto aH_2-bH_1$ and $(A_{1j},A_{2j})\mapsto aA_{2j}-bA_{1j}$ are surjective, with affine fibers. Normalizing the two original equations adds an invertible $2\times2$ coefficient matrix. In type~II, the maps $(R,H)\mapsto R(1,s)-sH(1,s)$ and $(B_i,A_i)\mapsto B_i(1,s)-sA_i(1,s)$ are surjective onto polynomials of degrees at most $2m$ and $m$, respectively. Their kernels and the coefficients used to normalize $xw,yw$ are free parameters on a coefficient chart. Thus each full coefficient field is a purely transcendental extension of the effective field. The finite incidences above stay integral after that extension, and Hilbert irreducibility applies on the full rational coefficient chart, intersected with its nonempty quasismooth open. This proves both assertions for the original families.
\end{proof}

\begin{Corollary}\label{cor:KP-generic-nonrational}
For every $k\geq1$, the generic members of the first two Kim--Park families are nonrational over their coefficient fields. Over a Hilbertian field $K$, a Hilbert subset of quasismooth members of either family consists of unirational, non-$K$-rational surfaces. Their conic bundle models are standard, have no rational section, and have discriminant degree $2k+5$.
\end{Corollary}
\begin{proof}
Combine Proposition~\ref{prop:KP-generic-standard}, Corollary~\ref{cor:KP-individual-nonrational}, and Proposition~\ref{prop:KP-bisections}.
\end{proof}

\begin{proof}[Proof of Theorem~\ref{thm:main}]
The classification discussed in Section~\ref{sec:introduction} reduces the assertion to the four families in \eqref{eq:families}. Corollaries~\ref{cor:JK-unirational} and~\ref{cor:KP3-unirational} prove unirationality for the Johnson--Koll\'ar and third Kim--Park family. Propositions~\ref{prop:KP1-conic}, \ref{prop:KP2-conic}, and~\ref{prop:KP-bisections} prove the conic bundle assertions and degree-two unirationality for the first two Kim--Park families. The explicit del Pezzo and cubic models are Propositions~\ref{prop:JK-model} and~\ref{prop:KP3-models}. Proposition~\ref{prop:JK-KP3-generic} and Corollary~\ref{cor:KP-generic-nonrational} prove generic nonrationality and its Hilbert specialisations. The relevant Hilbert subsets are Zariski dense in the quasismooth coefficient opens by Hilbert irreducibility.
\end{proof}

\subsection*{Positive characteristic}
The main theorem is stated in characteristic zero. The explicit unirationality proofs also work under the tame restrictions $\charac K\nmid2(2k+1)(4k+1)$ for the Johnson--Koll\'ar and second Kim--Park family, $\charac K\nmid2(k+1)(2k+1)$ for the first family, and $\charac K\nmid2(2k+1)(4k+1)(6k+1)$ for the third. The same local Jacobian checks and birational maps apply. 

For the conic models, a smooth total space rules out rank-one fibers and identically vanishing determinant by the local rank-two calculation in Lemma~\ref{lem:smooth-conic-discriminant}; the rational bisection remains separable since the characteristic is not two. For degree-two del Pezzo models over infinite fields, use Theorem~\ref{thm:Gao}; over finite fields, use \cite[Theorem~1.1]{FestiVanLuijk}. No positive-characteristic classification or generic nonrationality assertion is needed here.

\bibliographystyle{amsalpha}
\bibliography{Biblio}
\end{document}